\documentclass[prep,preprint]{imsart}

\RequirePackage{amsthm,amsmath,amsfonts,amssymb}
\RequirePackage[authoryear]{natbib}
\RequirePackage[colorlinks,citecolor=blue,urlcolor=blue]{hyperref}
\RequirePackage{graphicx}

\startlocaldefs
\usepackage{preamble}
\endlocaldefs

\begin{document}

\begin{frontmatter}
\title{High-dimensional Berry--Esseen bound for $m$-dependent random samples}
\runtitle{High-dimensional Berry--Esseen bound for $m$-dependent random samples}

\begin{aug}
\author[A]{\fnms{Heejong} \snm{Bong} \ead[label=e1]{hbong@andrew.cmu.edu}},
\author[A]{\fnms{Arun Kumar} \snm{Kuchibhotla} \ead[label=e2]{arunku@stat.cmu.edu}}
\and
\author[A]{\fnms{Alessandro} \snm{Rinaldo} \ead[label=e3,mark]{\{hbong, arunku, arinaldo\}@cmu.edu}}
\address[A]{Department of Statistics and Data Science,
Carnegie Mellon University,
\printead{e3}
}

\end{aug}

\begin{abstract}
In this work, we provide a $(n/m)^{-1/2}$-rate finite sample Berry--Esseen bound for $m$-dependent high-dimensional random vectors over the class of hyper-rectangles. This bound imposes minimal assumptions on the random vectors such as nondegenerate covariances and finite third moments. The proof uses inductive relationships between anti-concentration inequalities and Berry--Esseen bounds, which are inspired by the telescoping method of \citet{chen2004normal} and the recursion method of \citet{kuchibhotla2020high}. Performing a dual induction based on the relationships, we obtain tight Berry-Esseen bounds for dependent samples.
\end{abstract}

\begin{keyword}
\kwd{High-dimensional inference}
\kwd{Central limit theorem}
\kwd{Berry--Esseen bound}
\kwd{$m$-dependence structure}
\end{keyword}

\end{frontmatter}

\section{Introduction}

Recent advances in technology have led to the unprecedented availability of large-scale spatiotemporal data. An important challenge in the analyses of such data is to provide a theoretical guarantee of statistical inferences under temporal dependence. Many existing theoretical studies, such as \citet{liu2020matrix}, relied on parametric or distributional assumptions to give a valid confidence interval, but the validity of the assumptions remains questionable in real-world applications.  

In this paper, we focus on high-dimensional Central Limit Theorems (CLTs) over the class of hyper-rectangles. That is, our objective is to bound the Kolmogorov-Smirnov distance between a summation of samples $X_1, \dots, X_n \in \reals^p$ and its Gaussian approximation $Y$, denoted by
\begin{equation*}
    \mu\left( \sum_{i=1}^n X_i, Y \right) 
    := \sup_{r \in \reals^p} \abs{\Pr[\sum_{i=1}^n X_i \preceq r] - \Pr[Y \preceq r]},
\end{equation*}
where for two vectors $a, b \in \reals^p$, $a \preceq b$ means that $a_k \leq b_k$ for every $k \in [p]$. For independent samples, there has been a flurry of novel results since the seminal work of \citet{chernozhukov2013gaussian}. A popular approach has been the Lindeberg interpolation, leading to an $n^{-1/6}$ rate \citep{bentkus2000accuracy,lopes2022central}. Recently, \citet{kuchibhotla2020high} used a recursion method over the Lindeberg interpolation to establish a high-dimensional Berry--Esseen bound with the rate at most $(\log n)^{3/2} /\sqrt{n}$ under minimal assumptions: nondegenerate covariances and finite third moments of the samples. Here in this discussion, we are focusing on the dependence of the rate on the number $n$ of observations.
Our work extends their arguments to $m$-dependent cases in which $X_i \indep X_j$ if $\abs{i-j} > m$.

For $m$-dependent random samples, the optimal $(n/m)^{-1/2}$-rate has been only obtained in cases for $p=1$ by \citet{shergin1980convergence}. Later, the result was generalized to broad dependence structures by \citet{chen2004normal}. For high-dimensional cases, \citet{zhang2017gaussian} studied the Gaussian approximation for $m$-dependent sequences (see Appendix A.2 therein), but the rate of their Berry--Esseen bound and its dependency to sample assumptions were presented implicitly. More recently, \citet{chang2021central} provided a bound with the rate $n^{-1/6}$ for sub-exponential random vectors, using a decomposition of samples into ``large'' blocks and ``small'' blocks as done by \citet{romano2000more} for $1$-dimensional cases. 

We obtain a $(n/m)^{-1/2}$-rate bound under minimal assumptions that the samples have nondegenerate covariances and finite third moments. The proof uses inductive relationships between anti-concentration inequalities and Berry--Esseen bounds, which are results of the telescoping method by \citet{chen2004normal} and the recursion method by \citet{kuchibhotla2020high}. Performing a dual induction based on the relationships, we obtain tight Berry-Esseen bounds for dependent samples.

\section{Summation of $m$-dependent random vectors}

Let $X_1, X_2, \dots, X_n$ be $m$-dependent centered random vectors in $\reals^p$. That is, $X_i \indep X_j$ if $\abs{i-j} > m$, and $\Exp[X_i] = 0$ for $i \in [n]$. We use the superscript $X_i^{(k)}$ to notate each $k$-th element of $p$ dimensional random vector $X_i$. 
Let $Y_1, Y_2, \dots, Y_n$ be Gaussian random vectors with mean zero and the same second moment as $X_1, \dots, X_n$. i.e.,
\begin{equation*}
    \Exp[Y_i Y_j^\top] = \Exp[X_i X_j^\top],\quad\mbox{for all}\quad  i, j \in [n].
\end{equation*}
For the brevity of the following arguments, let
\begin{equation*}
    X_{[i,j)} := \{ X_k: i \leq k < j \}
    \textand 
    Y_{(i,j]} := \{ Y_k :i < k \leq j \}
\end{equation*}
and $S$ be the summation operator such that 
\begin{equation*}
    SX_{[i,j)} := \sum_{k: i \leq k < j} X_k
    \textand 
    SY_{(i,j]} := \sum_{k: i < k \leq j} Y_k.
\end{equation*}
For $i,j : 1 \leq i \leq j \leq n$, let ${\sigma}_{\min,[i,j]}$, $\underline{\sigma}_{[i,j]}$ and $\overline{\sigma}_{[i,j]}$ be
\begin{equation*} \begin{aligned}
    {\sigma}^2_{\min,[i,j]} :=
    & \min_{k \in [p]} \Var[SX^{(k)}_{[i,j]}], \\
    \underline{\sigma}^2_{[i,j]} :=
    & \lambda_{\min}(\Var[SX_{[i,j]}]), \\
    \overline{\sigma}^2_{[i,j]} :=
    & \lambda_{\max}(\Var[SX_{[i,j]}]).
\end{aligned} \end{equation*}
Also, let $\nu_1$ and $\nu_{3}$ be
\begin{equation*}
    \nu_{1} :=
    \mathrm{max} \left\{\begin{aligned}
    & \Exp[\norm{X_{i}}_\infty], \Exp[\norm{Y_{i}}_\infty]
    \end{aligned}:
    i \in [n] \right\},
\end{equation*}
\begin{equation*}
    \nu_{3} :=
    \mathrm{max} \left\{\begin{aligned}
    & \Exp[\norm{X_{i}}_\infty^3], \Exp[\norm{Y_{i}}_\infty^3]
    \end{aligned}:
    i \in [n] \right\}.
\end{equation*}

\vspace{0.1in}
{\bf Notation. } 
In the following argument, $C(\dots)$ is a constant with implicit dependency on the parameters in the parentheses, whose value is changing across lines. For absolute constants with no dependency, we omit the parentheses and denote them by $C$. $\mathbf{1}$ stands for the vector with elements $1$ in the appropriate dimension at each line. $\Ind\{\cdot\}$ denotes the indicator function.

\section{High-dimensional Berry--Esseen bound for $m$-dependent random samples with nondegenerate covariance matrices}

In this section, we state a high-dimensional Berry--Esseen bound for $m$-dependent random samples with nondegenerate covariance matrices. That is, the minimum and maximum eigenvalues of sums of consecutive samples are bounded away from $0$ and $\infty$, respectively. We assume 
\begin{align}
    \sigma^2_{\min,[i,j]} & \geq \sigma^2_{\min} \cdot (i-j+1) \cdot \min\{m, i-j+1\}, 
    \label{assmp:min_var} \tag{MIN-VAR} \\
    \underline{\sigma}^2_{[i,j]} & \geq \underline{\sigma}^2 \cdot (i-j+1) \cdot \min\{m, i-j+1\}, \label{assmp:min_ev} \tag{MIN-EV}\\ 
    \overline{\sigma}^2_{[i,j]} & \leq \overline{\sigma}^2 \cdot (i-j+1) \cdot \min\{m, i-j+1\}, 
    \label{assmp:max_ev} \tag{MAX-EV}
\end{align}
for arbitrary constants $\sigma_{\min}, \underline{\sigma}, \overline{\sigma} > 0$. 

\begin{theorem} \label{thm:1_dep_berry_esseen}
For $m=1$, suppose that Assumptions \eqref{assmp:min_var}, \eqref{assmp:min_ev} and \eqref{assmp:max_ev} hold. Then,
\begin{equation*}
\begin{aligned}
    & \mu\left( SX_{[1,n]}, SY_{[1,n]} \right) \\
    & \leq \frac{C}{\sqrt{n}} \left( 
        \frac{\nu_3}{\underline{\sigma}^3} \log(en) (\log(ep))^{2} \sqrt{\log(epn)}
        +  \frac{\nu_3^{1/3}\overline{\sigma}^{1/3}}{\sigma_{\min}\underline{\sigma}^{1/3}} \log(ep) \sqrt{\log(epn)}
    \right), \\
\end{aligned}
\end{equation*}
for some constant $C > 0$.
\end{theorem}

The Berry-Esseen bounds for $m>1$ come as corollaries of \cref{thm:1_dep_berry_esseen}, using a similar argument with Theorem 2 in \citet{shergin1980convergence} and Theorem 2.6 in \citet{chen2004normal}.

\begin{corollary} \label{thm:m_dep_berry_esseen}
Under the same assumptions as \cref{thm:1_dep_berry_esseen} except that $m = 1$ is replaced with $m\ge1$, with $n_{\mathrm{eff}} = n/m$, there exists a universal constant $C > 0$ such that
\begin{equation*}
\begin{aligned}
    & \mu\left( SX_{[1,n]}, SY_{[1,n]} \right) \\
    & \leq \frac{C}{\sqrt{n_\mathrm{eff}}} \left( 
        \frac{\nu_3}{\underline{\sigma}^3} \log(en_\mathrm{eff}) (\log(ep))^{2} \sqrt{\log(epn_\mathrm{eff})}
        +  \frac{\nu_3^{1/3}\overline{\sigma}^{1/3}}{\sigma_{\min}\underline{\sigma}^{1/3}} \log(ep) \sqrt{\log(epn_\mathrm{eff})}
    \right). \\
\end{aligned}
\end{equation*}
\end{corollary}


\section{Sketch of Proof} \label{sec:pf_sketch}

Here we provide a sketch of proof for the cases of $m=1$. We use inductive relationships between anti-concentration inequalities and Berry-Esseen bounds to prove \cref{thm:1_dep_berry_esseen}. For a random variable $X \in \reals^p$, anti-concentration inequalities give upperbounds for the probability of $X$ to be contained in a specific type of subsets in $\reals^p$. 
For example, \citet{nazarov2003maximal} and \citet{chernozhukov2017detailed} provided an anti-concentration inequality for Gaussian random variables to be contained in $A_{r,\delta}$, where for $r \in \reals^p$ and $\delta \in [0, \infty)$,
\begin{equation}
    A_{r,\delta} = \{x \in \reals^p: x \preceq r + \delta\mathbf{1}\}
    \setminus \{x \in \reals^p: x \preceq r - \delta\mathbf{1}\}.
\end{equation}

\begin{lemma}[Gaussian anti-concentration inequality; \citealp{nazarov2003maximal,chernozhukov2017detailed}]
\label{thm:G_anti_conc}
For a random vector $Y \dist N(0,\Sigma)$ in $\reals^p$, $r \in \reals^p$, and $\delta \in [0, \infty)$,
\begin{equation}
    \Pr[Y \in A_{r,\delta}] \leq C \delta \sqrt{\frac{\log(ep)}{\min_{i = 1, \dots, p} \Sigma_{ii}}}
\end{equation}
for an absolute constant $C > 0$.
\end{lemma}

We note that the anti-concentration inequality for $SX_{[1,n]}$ to be contained in $A_{r, \delta}$, i.e.,
\begin{equation} \label{eq:marginal_anti_conc}
    \Pr[SX_{[1,n]} \in A_{r,\delta}] \leq f(n, \delta),
\end{equation}
is a corollary of \cref{thm:G_anti_conc} and \cref{thm:1_dep_berry_esseen}. In cases of independent samples, \citet{kuchibhotla2020high} used the anti-concentration inequalities for $SX_{[1,i]}$ at $i < n$ in the Lindeberg swapping to prove the Berry-Esseen bounds for $SX_{[1,n]}$. In their proof, the dual induction between anti-concentration inequalities and Berry-Esseen bounds implied the desired result. However, the anti-concentration inequality for the marginal probability in \cref{eq:marginal_anti_conc} is not sufficient to prove the Berry--Esseen bound for dependent random vectors. Rather, we need a generalized anti-concentration inequality spanning the conditional probability,
\begin{equation}
    \Pr[SX_{[1,i]} \in A_{r,\delta} | X_{(i, n]}], ~~ \forall i \in [1,n],
\end{equation}
where $r, \delta$ are Borel-measurable functions with respect to $X_{(i, n]}$. Let
\begin{equation}
\begin{aligned}
    \kappa_i(\delta) := \sup_{r\in\reals^p} \Pr[SX_{[1,i]} \in A_{r,\delta} | X_{(i,n]}] \textand 
    \mu_i := \mu(SX_{[1,i]}, SY_{[1,i]}).
\end{aligned}
\end{equation}
The Lindeberg swapping method as in \citet{kuchibhotla2020high} derives the inductive relationship from $\kappa_i(\delta)$ to $\mu_i$: for $i \in [n]$ and $\vareps > 0$,
\begin{equation*}
    \sqrt{i} \mu_i
    \leq C(\nu_3, \sigma_{\min}, \underline{\sigma}) \vareps \log(ep)
    + C(\nu_3, \overline{\sigma}, \underline{\sigma}) \frac{ \log(en) (\log(ep))^{3/2}}{\vareps} \sup_{j: \frac{i}{2} < j < i} \sqrt{j} \kappa_{j}\left(\vareps\right).
\end{equation*}
On the other hand, the telescoping scheme in \citet{chen2004normal} derives the inductive relationship from $\mu_i$ to $\kappa_i(\delta)$: for $i \in [n]$,
\begin{equation*}
    \sqrt{i}\kappa_i(\delta) \leq C \left( \frac{\delta}{\sigma_{\min}} \sqrt{\log(ep)} + \max_{j: j \leq i-2} \sqrt{j} \mu_j \right).
\end{equation*}
The full statements of the inductive relationships are in \cref{thm:induction_lemma_1,thm:induction_lemma_2} of \cref{sec:pf_1_dep}. Summing the two equations and taking an appropriate $\vareps$, we obtain the following inductive inequality among $\mu_i$'s: for $i \geq 4$,
\begin{equation*}
    \sqrt{i} \mu_i \leq C(\nu_3, \sigma_{\min}, \overline{\sigma}, \underline{\sigma}) \log(en) (\log(ep))^{2} \sqrt{\log(epn)}, + \frac{1}{2} \max_{j: j \leq i-3} \sqrt{j} \mu_j.
\end{equation*}
Based on mathematical induction, we obtain for all $i\in[n]$,
\begin{equation*}
    \sqrt{i} \mu_i \leq 2 C(\nu_3, \sigma_{\min}, \overline{\sigma}, \underline{\sigma}) \log(en) (\log(ep))^{2} \sqrt{\log(epn)},
\end{equation*}
which proves the desired theorem.

\section{Discussion}

We derived a $n_\mathrm{eff}^{-1/2}$ scaling of Berry--Esseen bound for high-dimensional $m$-dependent random vectors over hyper-rectangles, where $n_\mathrm{eff} = n/m$ indicated the effective sample size. This result only required nondegenerate covariances and finite third moments of the random vectors. The exact rate with respect to $n_\mathrm{eff}$ was $(\log n_\mathrm{eff})^{-3/2}/\sqrt{n_\mathrm{eff}}$, where the dependency on the dimension $p$ was logarithmic. The rate had the same scale up to logarithm as the result of \citet{shergin1980convergence} on $1$-dimensional $m$-dependent samples, which is also known to be unimprovable based on \citet{berk1973central}. Our result supports the high-dimensional CLT and use of bootstrap over hyper-rectangles under $m$-dependency between samples.

Our advancement in the Gaussian approximation rate of $m$-dependent samples could benefit the theoretical analyses under physical dependence frameworks. \citet{zhang2018gaussian} introduced the $m$-approximation technique to study the Gaussian approximation of weakly dependent time series under physical dependence. The technique extends the Berry--Esseen bounds for $m$-dependent samples to weaker temporal dependencies (see Theorem 2.1 and the end of Section 2.2 therein). Similarly, \citet{chang2021central} extended the $n^{-1/6}$ rate under $m$-dependence to samples with physical dependence. The resulting rate in Theorem 3 was better than the best rate of \citet{zhang2017gaussian} at that time.

Another important future direction is extending our technique to samples with generalized graph dependency. Random vectors $X_1, X_2, \dots, X_n \in \reals^p$ are said to have dependency structure defined by graph $G = ([n], E)$ if $X_i \indep X_j$ if $(i,j) \in E$. Graph dependency generalizes $m$-dependence as a special case by taking $E = \{(i,j): \abs{i-j} \leq m\}$. The only CLT result up to our best knowledge has been \citet{chen2004normal} for $1$-dimensional samples with graph dependency. Extending their result to high-dimensional samples has a huge potential to advance statistical analyses on network data, which is another data type with increasing availability.

\bibliographystyle{imsart-nameyear} 
\bibliography{aoas-refs}       





\newpage
\begin{appendix}

\section{Proofs}

\subsection{Detailed proof of \cref{thm:1_dep_berry_esseen}} \label{sec:pf_1_dep}

Let our induction hypothesis be
\begin{equation*}
\begin{aligned}
    \sqrt{i}\mu_i
    & \leq C \left( 
        \frac{\nu_3}{\underline{\sigma}^3} \log(en) \log(ep)
        +  \frac{\nu_1\overline{\sigma}^{1/3}}{\sigma_{\min}\underline{\sigma}^{1/3}} 
    \right) 
    \log(ep) \sqrt{\log(epn)},
\end{aligned}
\end{equation*}
for $i\in[n]$ and some absolute constant $C > 0$. Since $\nu_1 \geq \sigma_{\min} \geq \underline{\sigma} > 0$, every term inside the parentheses on the right hand side is larger than $1$. Thus, if $C \geq 1$, the induction hypothesis, requiring $\mu_i \leq 1$, trivially holds for $i = 1, 2, 3$. Now we apply induction for $i \geq 4$ using the induction lemmas in \cref{sec:pf_sketch}. Here are the full statements of the lemmas.

\begin{lemma} 
\label{thm:induction_lemma_1}
    If Assumptions \eqref{assmp:min_var}, \eqref{assmp:min_ev} and \eqref{assmp:max_ev} with $m=1$ hold, then for $i \in [n]$,
    \begin{equation} \label{eq:induction_1}
        \kappa_i(\delta) \leq \frac{C}{\sqrt{\max\{i-2,1\}}} \left( \frac{\delta + \nu_1}{\sigma_{\min}} \sqrt{\log(ep)} + \max_{j: j \leq i-2} \sqrt{j} \mu_j \right).
    \end{equation}
\end{lemma}

\begin{lemma} 
\label{thm:induction_lemma_2}
    If Assumptions \eqref{assmp:min_var}, \eqref{assmp:min_ev} and \eqref{assmp:max_ev} with $m=1$ hold, then for any $\vareps > 0$,
    \begin{equation} \label{eq:induction_2}
    \begin{aligned}
        \mu_n
        & \leq\frac{C}{\sqrt{n}} \left( \frac{\nu_3 (\log(ep))^{3/2}}{\underline{\sigma}^3} 
        + \frac{\vareps \log(ep)}{\sigma_{\min}} \right)
        + \frac{C}{pn} \nu_3 \left(\frac{\overline{\sigma}}{\vareps^3 \underline{\sigma}} + \frac{1}{\vareps\underline\sigma^2}\right) \\
        & \quad + C \nu_3  \left( 
            \frac{\overline{\sigma}}{\vareps^2 \underline{\sigma}} + \frac{1}{\underline\sigma^2} \log\left(1 + \frac{\sqrt{n}\underline{\sigma}}{\vareps}\right)
        \right) (\log(ep))^{3/2}
        \sup_{i: \frac{n}{2}+1 < i \leq n} \frac{\kappa_{i-1}(\delta_i)}{\vareps_i},
    \end{aligned}
    \end{equation}
    where $\vareps_i^2 = \vareps^2 + \max\{\underline{\sigma}^2 (n-i) - \overline{\sigma}^2, 0\}$ and $\delta_i = C \vareps_i \sqrt{\log(pn)}$.
\end{lemma}

Taking \cref{eq:induction_1,eq:induction_2} together,
\begin{equation*}
\begin{aligned}
    \sqrt{n}\mu_n
    & \leq C \left( 
        \frac{\nu_3 (\log(ep))^{3/2}}{\underline{\sigma}^3} 
        + \frac{\vareps \log(ep)}{\sigma_{\min}} 
    \right) 
    + \frac{C}{p\sqrt{n}} \nu_3 \left(
        \frac{\overline{\sigma}}{\vareps^3 \underline{\sigma}} 
        + \frac{1}{\vareps\underline\sigma^2}
    \right)\\
    & \quad + C \nu_3 \left( 
            \frac{\overline{\sigma}}{\vareps^2 \underline{\sigma}} 
            + \frac{1}{\underline\sigma^2} \log\left(1 + \frac{\sqrt{n}\underline{\sigma}}{\vareps}\right)
    \right) (\log(ep))^{3/2}
    \sup_{i: \frac{n}{2}+1 < i \leq n} \frac{\sqrt{n} \kappa_{i-1}(\delta_i)}{\vareps_i} \\
    & \leq C \left( 
        \frac{\nu_3 (\log(ep))^{3/2}}{\underline{\sigma}^3} 
        + \frac{\vareps \log(ep)}{\sigma_{\min}} 
    \right) 
    + \frac{C}{p\sqrt{n}} \nu_3 \left(
        \frac{\overline{\sigma}}{\vareps^3 \underline{\sigma}} 
        + \frac{1}{\vareps\underline\sigma^2}
    \right)\\
    & \quad + C \nu_3 \left( 
            \frac{\overline{\sigma}}{\vareps^2 \underline{\sigma}} 
            + \frac{1}{\underline\sigma^2} \log\left(1 + \frac{\sqrt{n}\underline{\sigma}}{\vareps}\right)
    \right) (\log(ep))^{3/2} \\
    & \quad \times \sup_{i: \frac{n}{2}+1 < i \leq n} \frac{1}{\vareps_i} \left( 
        \frac{C \vareps_i \sqrt{\log(pn)} + \nu_1}{\sigma_{\min}} \sqrt{\log(ep)} 
        + \max_{j: j \leq i-3} \sqrt{j} \mu_j 
    \right) \\
    & \leq C \left( 
        \frac{\nu_3 (\log(ep))^{3/2}}{\underline{\sigma}^3} 
        + \frac{\vareps \log(ep)}{\sigma_{\min}} 
    \right) 
    + \frac{C}{p\sqrt{n}} \nu_3 \left(
        \frac{\overline{\sigma}}{\vareps^3 \underline{\sigma}} 
        + \frac{1}{\vareps\underline\sigma^2}
    \right)\\
    & \quad + C \frac{\nu_3}{\sigma_{\min}} \left( 
            \frac{\overline{\sigma}}{\vareps^2 \underline{\sigma}} 
            + \frac{1}{\underline\sigma^2} \log\left(1 + \frac{\sqrt{n}\underline{\sigma}}{\vareps}\right)
    \right)
        \sqrt{\log(pn)} (\log(ep))^2 \\
    & \quad + C \frac{\nu_3}{\vareps} \left( 
            \frac{\overline{\sigma}}{\vareps^2 \underline{\sigma}} 
            + \frac{1}{\underline\sigma^2} \log\left(1 + \frac{\sqrt{n}\underline{\sigma}}{\vareps}\right)
    \right) (\log(ep))^{3/2} 
    \max_{j: j \leq n-3} \sqrt{j} \mu_j, \\
\end{aligned}
\end{equation*}
where $\vareps \geq \nu_1$. We can take $\vareps = C \left( (\nu_3/\underline{\sigma}^2) \log(en)(\log(ep))^{3/2} + (\nu_3 \overline{\sigma}/\underline{\sigma})^{1/3} \sqrt{\log(ep)} \right)$ for some $C \geq 1$ such that
\begin{equation*}
\begin{aligned}
    \sqrt{n}\mu_n
    & \leq C \frac{\nu_3}{\underline{\sigma}^3} \log(en) (\log(ep))^{5/2}
    + C \frac{\nu_3^{1/3}\overline{\sigma}^{1/3}}{\sigma_{\min}\underline{\sigma}^{1/3}} (\log(ep))^{3/2}
    + \frac{C}{\sqrt{n}} \frac{1}{p\sqrt{\log(ep)}}\\
    & \quad + C \frac{\nu_3^{1/3}\overline{\sigma}^{1/3}}{\sigma_{\min}\underline{\sigma}^{1/3}} \sqrt{\log(pn)} \log(ep) 
    + C \frac{\nu_3}{\sigma_{\min}\underline{\sigma}^2} \log(en) \sqrt{\log(pn)} (\log(ep))^2\\
    & \quad + \frac{1}{2} \max_{j: j \leq n-3} \sqrt{j} \mu_j. \\
\end{aligned}
\end{equation*}
Then, by the induction hypothesis,
\begin{equation*}
\begin{aligned}
    \sqrt{n}\mu_n
    & \leq C \left( 
        \frac{\nu_3}{\underline{\sigma}^3} \log(en) (\log(ep))^{2} \sqrt{\log(epn)}
        +  \frac{\nu_3^{1/3}\overline{\sigma}^{1/3}}{\sigma_{\min}\underline{\sigma}^{1/3}} \log(ep) \sqrt{\log(epn)}
    \right), \\
\end{aligned}
\end{equation*}
which proves our theorem.

\subsection{Proof of \cref{thm:m_dep_berry_esseen}}

Let $n' = \floor{\frac{n-1}{m}}$. Suppose that $X'_i = \frac{1}{m} SX_{((i-1)m, im]}$ for $i \in [n'-1]$ and $X'_{n'} = \frac{1}{m} SX_{((n'-1)m, n]}$. We define $Y'_i$ similarly for $i \in [n']$. Then, $X'_1, \dots, X'_{n'}$ are $1$-dependent random vectors in $\reals^p$. The corollary follows the observations that $X'_i$ and $Y'_i$ satisfy Assumptions \eqref{assmp:min_var}, \eqref{assmp:min_ev} and \eqref{assmp:max_ev} with $m=1$ and the same $\sigma_{\min}$, $\underline{\sigma}$, $2\overline{\sigma}$, and 2$\nu_3$.

\subsection{Proof of \cref{thm:induction_lemma_1}}

If $i \leq 4$, because $\nu_1 \geq \sigma_{\min}$, the righthand side of \cref{eq:induction_1} is greater than $1$ for $C \geq 2$, where the lemma holds trivially.

For $i \geq 3$, we use the telescoping scheme in Proposition 3.2, \citet{chen2004normal} to prove the anti-concentration inequality. Let our induction hypothesis be
\begin{equation*}
    \sqrt{i-2} \kappa_i\left(2 \sqrt{\frac{i}{i-2}} \nu_1\right) 
    \leq C \left( \frac{\nu_1}{\sigma_{\min}} \sqrt{\log(ep)} + \sup_{j \leq i-2} \sqrt{j} \mu_{j} \right)
\end{equation*}
for some absolute constant $C > 0$. If $C \geq 2$, the righthand side is greater than $2$, so the induction hypothesis trivially holds for $i=3,4$. Suppose that $i \geq 5$. Conditional on $X_{(i,n]}$, $r$ and $\delta$ are almost surely constant. For $\vareps > 0$, let $f_\vareps$ be a function $\reals^p \rightarrow \reals$ such that
\begin{equation*}
    f_\vareps(x) = \begin{cases}
        \frac{w+\delta+\vareps}{\vareps}, & \text{for} ~~ - \delta - \vareps < w \leq - \delta, \\
        1, & \text{for} ~~ -\delta < w \leq \delta, \\
        \frac{\delta+\vareps-w}{\vareps}, & \text{for} ~~ \delta < w \leq \delta + \vareps, \\
        0, & \text{elsewhere,}
    \end{cases}
\end{equation*}
where $w = \max\{x^{(k)} - r^{(k)}: k \in [p]\}$. We note that 
\begin{equation} \label{eq:grad_f}
\begin{aligned}
    \norm{\nabla f_\vareps(x)}_1 
    & = \begin{cases}
        \frac{1}{\vareps}, 
        & \text{for} ~~ - \delta - \vareps < w \leq - \delta ~\text{or}~ \delta < w \leq \delta + \vareps, \\
        0, & \text{elsewhere}
    \end{cases} \\
    & = \frac{1}{\vareps} \left( 
        \Ind\{x \in A_{r-(\delta+\vareps/2)\mathbf{1}, \vareps/2}\}
        + \Ind\{x \in A_{r+(\delta+\vareps/2)\mathbf{1}, \vareps/2}\}
    \right),
\end{aligned}
\end{equation}
almost everywhere. Because $f_\vareps(x) \geq \Ind\{x \in A_{r,\delta}\}$,
\begin{equation*}
\begin{aligned}
    \Pr[SX_{[1,i]} \in A_{r,\delta} | X_{(i,n]}]
    & \leq \Exp[f_\vareps(SX_{[1,i]}) | X_{(i,n]}] \\
    & = \Exp[f_\vareps(SX_{[1,i]}) - f_\vareps(SX_{[1,i]} - X_{i-1}) | X_{(i,n]}] \\
    & \quad + \Exp[f_\vareps(SX_{[1,i]} - X_{i-1}) | X_{(i,n]}] \\
    & =: H_1 + H_2.
\end{aligned}
\end{equation*}
Because of Taylor's expansion and H\"older's inequality,
\begin{equation*}
\begin{aligned}
    H_1 
    & \leq \Exp\left[ \left. \int_0^1 (1-t) \inner{\nabla f_\vareps(SX_{[1,i]} - t X_{i-1}), X_{i-1}} dt \right| X_{(i,n]} \right] \\
    & \leq \Exp\left[ \left. \int_0^1 \norm{\nabla f_\vareps(SX_{[1,i]} - t X_{i-1})}_1 \norm{X_{i-1}}_\infty dt \right| X_{(i,n]} \right] \\
    & \leq \frac{1}{\vareps} \Exp\left[ \left. \int_0^1 \left(
    \begin{aligned}
        & \Pr[SX_{[1,i-2]} \in A_{r_1, \vareps/2} | X_{(i-2,n]}] \\
        & + \Pr[SX_{[1,i-2]} \in A_{r_2, \vareps/2} | X_{(i-2,n]}]
    \end{aligned} \right)
    \norm{X_{i-1}}_\infty dt \right| X_{(i,n]} \right] \\
    & \leq \frac{1}{\vareps} \Exp\left[ \left. \int_0^1 2 \kappa_{i-2}\left(\frac{\vareps}{2}\right) \norm{X_{i-1}}_\infty \right| X_{(i,n]} \right] \\
    & = \frac{2}{\vareps} \kappa_{i-2}\left(\frac{\vareps}{2}\right) \Exp[\norm{X_{i-1}}_\infty]
    \leq \frac{2\nu_1}{\vareps} \kappa_{i-2}\left(\frac{\vareps}{2}\right)
\end{aligned}
\end{equation*}
where $r_1 = r - (1-t) X_{i-1} - X_i -(\delta+\vareps/2)\mathbf{1}$ and $r_2 = r - (1-t) X_{i-1} - X_i + (\delta+\vareps/2)\mathbf{1}$ are Borel measurable functions with respect to $X_{(i-2,n]}$. On the other hand,
\begin{equation*}
\begin{aligned}
    H_2
    & = \Exp[f_\vareps(SX_{[1,i-2]} + X_{i}) | X_{(i,n]}] \\
    & \leq \Exp[\Pr[SX_{[1,i-2]} \in A_{r_3, \delta} | X_{[i,n]}] | X_{(i,n]}]
\end{aligned}
\end{equation*}
where $r_3 = r - X_i$ is a Borel measurable function with respect to $X_{[i,n]}$. Because $X_{[1,i-2]} \indep X_{[i,n]}$, 
\begin{equation*}
\begin{aligned}
    \Pr[SX_{[1,i-2]} \in A_{r_3, \delta} | X_{[i,n]}] 
    & \leq \Pr[SY_{[1,i-2]} \in A_{r_3, \delta} | X_{[i,n]}] 
    + \mu(SX_{[1,i-2]}, SY_{[1,i-2]}) \\
    & \leq \frac{C \delta}{\sigma_{\min}} \sqrt{\frac{\log(ep)}{i-2}} + \mu_{i-2}.
\end{aligned}
\end{equation*}
almost surely due to the Gaussian anti-concentration inequality (\cref{thm:G_anti_conc}). In sum,
\begin{equation*}
    \Pr[SX_{[1,i]} \in A_{r,\delta} | X_{(i,n]}]
    \leq \frac{2\nu_1}{\vareps} \kappa_{i-2}\left(\frac{\vareps}{2}\right)
    + \frac{C \delta}{\sigma_{\min}} \sqrt{\frac{\log(ep)}{i-2}} + \mu_{i-2}.
\end{equation*}
Because the righthand side is not dependent on $r$, for any $\vareps, \delta > 0$, 
\begin{equation} \label{eq:K_i_delta}
    \kappa_i(\delta) \leq \frac{2\nu_1}{\vareps} \kappa_{i-2}\left(\frac{\vareps}{2}\right)
    + \frac{C \delta}{\sigma_{\min}} \sqrt{\frac{\log(ep)}{i-2}} + \mu_{i-2},
\end{equation}
where $C > 0$ is an absolute constant. Taking $\vareps = 4\sqrt{\frac{i-2}{i-4}} \nu_1$ and $\delta = 2 \sqrt{\frac{i}{i-2}} \nu_1$,
\begin{equation*}
\begin{aligned}
    \sqrt{i-2} \kappa_i\left(2 \sqrt{\frac{i}{i-2}} \nu_1\right) 
    & \leq \frac{1}{2} \sqrt{i-4} \kappa_{i-2}\left(2 \sqrt{\frac{i-2}{i-4}} \nu_1\right) \\
    & \quad + \frac{C \nu_1}{\sigma_{\min}} \sqrt{\frac{i \cdot \log(ep)}{i-2}} + \sqrt{i-2} \mu_{i-2}.
\end{aligned}
\end{equation*}
By the induction hypothesis at $i-2 \geq 3$, 
\begin{equation*}
\begin{aligned}
    \sqrt{i-2} \kappa_i\left(2 \sqrt{\frac{i}{i-2}} \nu_1\right) 
    & \leq  C \left(\frac{ \nu_1}{\sigma_{\min}} \sqrt{\log(ep)} + \sup_{j \leq i-2} \sqrt{j} \mu_{j} \right),
\end{aligned}
\end{equation*}
which validates the hypothesis at $i$.

Furthermore, for any $i \geq 5$ and $\delta > 0$, taking $\vareps = 4\sqrt{\frac{i-2}{i-4}} \nu_1$ in \cref{eq:K_i_delta},
\begin{equation*}
\begin{aligned}
    \kappa_i(\delta) 
    & \leq \frac{1}{2} \sqrt{\frac{i-4}{i-2}} \kappa_{i-2}\left(2 \sqrt{\frac{i-2}{i-4}} \nu_1\right)
    + \frac{C \delta}{\sigma_{\min}} \sqrt{\frac{\log(ep)}{i-2}} + \mu_{i-2} \\
    & \leq \frac{C}{i-2} \left(\frac{\nu_1 + \delta}{\sigma_{\min}} \sqrt{\log(ep)} + \sup_{j \leq i-2} \sqrt{j} \mu_{j} \right),
\end{aligned}
\end{equation*}
which proves the lemma.

\subsection{Proof of \cref{thm:induction_lemma_2}}
For Berry-Esseen bounds, we use the following smoothing lemma.
\begin{lemma}[Lemma 1, \citealp{kuchibhotla2020high}]
\label{thm:G_smoothing}
Suppose that $X$ is a $p$-dimensional random vector, and $Y \dist N(0,\Sigma)$ is a $p$-dimensional Gaussian random vector. Then, for any $\vareps > 0$ and a standard Gaussian random vector $Z$, 
\begin{equation}
    \mu(X, Y) \leq C \mu(X + \varepsilon Z, Y + \varepsilon Z) + \frac{C}{\sqrt{\min_{i = 1, \dots, p} \Sigma_{ii}}} \vareps \log(ep).
\end{equation}
\end{lemma}

Let  
\begin{equation*}
    W_{[i,j]} := X_{[1,i)} \cup Y_{(j,n]},
    W_{[i,j]}^{\indep} := W_{[i-1,j+1]},
\end{equation*}
and $\varphi_{r,\vareps}(x) = \Exp[ x + \vareps Z \preceq r ]$ for $r \in \reals^p, \vareps > 0$, where $Z$ is the standard Gaussian random variable in $\reals^p$.
For any $\vareps > 0$, by Lindeberg's swapping,
\begin{equation} \label{eq:Lindeberg_swap}
\begin{aligned}
    & \mu\left(SX_{[1,n]} + \vareps Z, SY_{[1,n]} + \vareps Z \right) \\
    & = \sup_{r \in \reals^p} \abs*{\varphi_{r,\vareps}\left( SX_{[1,n]} \right) 
    - \varphi_{r,\vareps}(SY_{[1,n]})} \\
    & = \sup_{r \in \reals^p} \abs*{ \sum_{j=1}^n 
        \Exp \left[
        \varphi_{r,\vareps}(SW_{[j,j]} + X_j)
        - \varphi_{r,\vareps}(SW_{[j,j]} + Y_j)
        \right]}.
\end{aligned}
\end{equation}
By Taylor's expansion, for $j = 1, \dots, n$,
\begin{equation*}
\begin{aligned}
    & \varphi_{r,\vareps}(SW_{[j,j]} + X_j)\\
    & = \varphi_{r,\vareps}(SW_{[j,j]})
    + \inner*{ \nabla \varphi_{r,\vareps}(SW_{[j,j]}), X_j}
    + \frac{1}{2} \inner*{ 
    \nabla^2 \varphi_{r,\vareps}(SW_{[j,j]}), X_j^{\otimes 2}} \\
    & \quad + \frac{1}{6} \int_0^1 (1-t)^3 \inner*{
    \nabla^3 \varphi_{r,\vareps}(SW_{[j,j]} + t X_j), X_j^{\otimes 3}} ~dt.
\end{aligned}
\end{equation*}
We further apply Taylor's expansion to the second and third terms:
\begin{equation*}
\begin{aligned}
    & \inner*{ \nabla \varphi_{r,\vareps}(SW_{[j,j]}), X_j} \\
    & = \inner*{ \nabla \varphi_{r,\vareps}(SW_{[j,j]}^{\indep}), X_j} 
    + \inner*{ 
    \nabla^2 \varphi_{r,\vareps}(SW_{[j,j]}^{\indep}),   X_j 
    \otimes (X_{j-1} + Y_{j+1})} \\
    & \quad + \frac{1}{2} \int_0^1 (1-t)^2 \left\langle
        \nabla^3 \varphi_{r,\vareps}\left(
            SW_{[j,j]}^{\indep} + t (X_{j-1} + Y_{j+1}) 
        \right), 
        X_j \otimes (X_{j-1} + Y_{j+1})^{\otimes 2}
    \right\rangle ~dt,
\end{aligned}
\end{equation*}
\begin{equation*}
\begin{aligned}
    & \inner*{ 
        \nabla^2 \varphi_{r,\vareps}(SW_{[j,j]}), X_j^{\otimes 2}} \\
    & = \inner*{\nabla^2 \varphi_{r,\vareps}(SW_{[j,j]}^{\indep}), X_j^{\otimes 2}} \\
    & + \int_0^1 (1-t) \left\langle
    \nabla^3 \varphi_{r,\vareps}\left(
          SW_{[j,j]}^{\indep} + t (X_{j-1} + Y_{j+1})
    \right), 
    X_j^{\otimes 2} \otimes (X_{j-1} + Y_{j+1})\right\rangle ~dt.
\end{aligned}
\end{equation*}

Last, 
\begin{equation*}
\begin{aligned}
    & \inner*{ 
        \nabla^2 \varphi_{r,\vareps}(SW_{[j,j]}^{\indep}), 
        X_j \otimes (X_{j-1}+Y_{j+1})} \\
    & = \inner*{ 
        \nabla^2 \varphi_{r,\vareps}(SW_{[j,j]}^{\indep}),   
        X_j \otimes X_{j-1}} 
    + \inner*{ 
        \nabla^2 \varphi_{r,\vareps}(SW_{[j,j]}^{\indep}), 
        X_j \otimes Y_{j+1}} \\
    & = \inner*{ 
        \nabla^2 \varphi_{r,\vareps}(SW_{[j-1,j]}^{\indep}), 
        X_j \otimes X_{j-1}}
    + \inner*{ 
        \nabla^2 \varphi_{r,\vareps}(SW_{[j+1,j]}^{\indep}),   
        X_j \otimes Y_{j+1}} \\
    & + \int_0^1 (1-t) \inner*{
        \nabla^3 \varphi_{r,\vareps}(SW_{[j-1,j]}^{\indep} + t X_{j-2}), 
        X_j \otimes X_{j-1} \otimes X_{j-2}} ~dt \\
    & + \int_0^1 (1-t) \inner*{
    \nabla^3 \varphi_{r,\vareps}(SW_{[j+1,j]}^{\indep}+ t Y_{j+2}),
        X_j \otimes Y_{j+1} \otimes Y_{j+2})} ~dt.
\end{aligned}
\end{equation*}
In sum,
\begin{equation*}
\begin{aligned}
    & \Exp \left[
        \varphi_{r,\vareps}(SW_{[j,j]} + X_j)\right] \\
    & = \Exp \left[
    \begin{aligned}
        & \varphi_{r,\vareps}(SW_{[j,j]})
        + \inner*{ \nabla \varphi_{r,\vareps}(SW_{[j,j]}^{\indep}), X_j } 
        + \frac{1}{2} 
        \inner*{\nabla^2 \varphi_{r,\vareps}(SW_{[j,j]}^{\indep}), X_j^{\otimes 2}} \\
        & + \inner*{ 
            \nabla^2 \varphi_{r,\vareps}(SW_{[j-1,j]}^{\indep}),  
            X_j \otimes X_{j-1}} 
        + \inner*{ 
            \nabla^2 \varphi_{r,\vareps}(SW_{[j+1,j]}^{\indep}),
            X_j \otimes Y_{j+1}} 
        + \mathfrak{R}_3
    \end{aligned} \right] \\
    & = \Exp\left[ \varphi_{r,\vareps}(SW_{[j,j]}) \right]
    + \inner*{ \Exp\left[ \nabla \varphi_{r,\vareps}(SW_{[j,j]}^{\indep}) \right], 
        \Exp\left[ X_j \right]}
    +
        \frac{1}{2} 
        \inner*{ \Exp\left[ \nabla^2 \varphi_{r,\vareps}(SW_{[j,j]}^{\indep}) \right], 
        \Exp\left[ X_j^{\otimes 2}\right]} \\
    & \quad + \inner*{ 
        \Exp\left[ \nabla^2 \varphi_{r,\vareps}(SW_{[j-1,j]}^{\indep}) \right], 
        \Exp\left[ X_j \otimes X_{j-1} \right]}
    + \Exp[\mathfrak{R}_{X_j,3}]
\end{aligned}
\end{equation*}
where $\mathfrak{R}_{X_j,3}$ is the remainder term involving with the third moments of the random vectors. Due to the moment matching between $X_i$ and $Y_i$ up to the second order,
\begin{equation*}
    \sum_j \Exp \left[
        \varphi_{r,\vareps}(SW_{[j,j]} + X_j)
        - \varphi_{r,\vareps}(SW_{[j,j]} + Y_j) \right]
    = \sum_j \Exp \left[ \mathfrak{R}_{X_j, 3} 
        - \mathfrak{R}_{Y_j, 3} \right],
\end{equation*}
and
\begin{equation*}
    \mu\left( SX_{[1,n]} + \vareps Z,   SY_{[1,n]} + \vareps Z \right)
    \leq \sum_j \sup_{r \in \reals^p} \abs*{ \Exp \left[ \mathfrak{R}_{X_j, 3} \right] }
    + \sum_j \sup_{r \in \reals^p} \abs*{ \Exp \left[ \mathfrak{R}_{Y_j, 3} \right] }.
\end{equation*}

With $n_0=\floor{n/2}$,
\begin{equation*}
    \sum_j \sup_{r \in \reals^p} \abs*{ \Exp \left[ \mathfrak{R}_{X_j, 3} \right] }
    \leq \sum_{j=n-n_0+1}^n \sup_{r \in \reals^p} \abs*{ \Exp \left[ \mathfrak{R}_{X_j, 3} \right] }
    + \sum_{j=1}^{n-n_0} \sup_{r \in \reals^p} \abs*{ \Exp \left[ \mathfrak{R}_{X_j, 3} \right] }.
\end{equation*}
For the remainder terms with $X_j^{\otimes 3}$,
\begin{equation*}
\begin{aligned}
    & \abs*{ \Exp\left[ \int_0^1 (1-t)^3 \inner*{
    \nabla^3 \varphi_{r,\vareps}(SW_{[j,j]} + t X_j), X_j^{\otimes 3}} ~dt
    \right] }\\
    & =  \int_0^1 \abs*{ \Exp\left[  \inner*{
    \nabla^3 \varphi_{r,\vareps}(SW_{[j,j]} + t X_j), X_j^{\otimes 3}}
    \right] } ~dt\\
    & \leq \int_0^1 \Exp\left[  \norm*{
    \nabla^3 \varphi_{r,\vareps}(SW_{[j,j]} + t X_j)}_1 \norm*{X_j}_\infty^3
    \right]  ~dt.\\
\end{aligned}
\end{equation*}
We note that by the Schur complement,
\begin{equation*}
\begin{aligned}
    \Var[SY_{(j,n]}|Y_{j}] 
    & = \Var[SY_{(j,n]}] - \Cov[Y_{j+1},Y_j] \Var[Y_j]^{-1} \Cov[Y_j,Y_{j+1}]\\
    & = \Var[SY_{(j,n]}] - \Var[Y_{j+1}] + \Var[Y_{j+1}|Y_j] \\
    & \succeq \Var[SY_{(j+1,n]}] - \Var[Y_{j+1}] 
    \succeq \max\{\underline{\sigma}^2 (n-j) - \overline{\sigma}^2, 0\} I_p
\end{aligned}
\end{equation*}
and that 
\begin{equation*}
    SY_{(j,n]}|Y_{j} \overset{d}{=} S^o + \sqrt{\max\{\underline{\sigma}^2 (n-j) - \overline{\sigma}^2 , 0\}} \cdot Z
\end{equation*}
alsmot surely where $S^o$ is the normal random variable with mean $\Exp[SY_{(j,n]}|Y_j]$ and variance $\Var[SY_{(j,n]}|Y_j]- \max\{\underline{\sigma}^2 (n-j) - \overline{\sigma}^2 , 0\} I_d$.
For $\vareps_j^2 = \vareps^2 + \max\{\underline{\sigma}^2 (n-j) - \overline{\sigma}^2, 0\}$,
\begin{equation*}
\begin{aligned}
    & \Exp\left[ \norm*{\nabla^3 \varphi_{r,\vareps}(SW_{[j,j]} + t X_j)}_1 
    \norm*{X_j}_\infty^3 \right] 
    = \Exp\left[ \norm*{\nabla^3 \varphi_{r,\vareps_j}(SX_{[1,j)} + S^o + t X_j)}_1 
    \norm{X_j}_\infty^3 \right]. \\
\end{aligned}
\end{equation*}
For $j=1, \dots, n_0$,
\begin{equation*}
\begin{aligned}
    \Exp\left[ \norm*{\nabla^3 \varphi_{r, \vareps_j}(SX_{[1,j)} + S^o + t X_j)}_1 
    \norm{X_j}_\infty^3 \right]
    & \leq C \frac{(\log(ep))^{3/2}}{\vareps_j^3}
    \Exp[\norm{X_j}_\infty^3] \\
    & \leq C \frac{(\log(ep))^{3/2}}{\underline{\sigma}^3} \frac{\nu_3}{\sqrt{n}^3}.
\end{aligned}
\end{equation*}
For $j=n_0+1, \dots, n$, based on the proof of Lemma 2 in \citet{kuchibhotla2020high},
\begin{equation*}
\begin{aligned}
    & \Exp\left[ \norm*{\nabla^3 \varphi_{r,\vareps_j}(SX_{[1,j)} + S^o + t X_j)}_1 
    \norm{X_j}_\infty^3 \right] \\
    & \leq \Exp\left[ C \left( 
        \frac{1}{\vareps_j^3 p n}
        + \frac{(\log(ep))^{3/2}}{\vareps_j^3} 
        \Ind \left\{ SX_{[1,j)} + S^o + t X_j \in A_{r, \delta_j} \right\}
    \right)
    \cdot \norm{X_j}_\infty^3 \right]\\
    & \leq \Exp\left[ C \left(
        \frac{1}{\vareps_j^3 p n}
        + \frac{(\log(ep))^{3/2}}{\vareps_j^3} 
        \Pr[ SX_{[1,j)} \in A_{r - S^o - t X_j, \delta_j} | X_j, S^o]
    \right)
    \cdot \norm{X_j}_\infty^3 \right]\\
    & \leq C \nu_3 \left( 
        \frac{1}{\vareps_j^3 p n}
        + \frac{(\log(ep))^{3/2}}{\vareps_j^3} \kappa_{j-1}(\delta_j) 
    \right)\\
\end{aligned}
\end{equation*}
where $\delta_j = c \vareps_j \sqrt{\log(pn)}$. Applying similar arguments to the other remainder terms, \cref{thm:G_smoothing} obtains
\begin{equation*}
\begin{aligned}
    & \mu\left( SX_{[1,n]}, SY_{[1,n]} \right) \\
    & \leq \mu\left(SX_{[1,n]} + \vareps Z, SY_{[1,n]} + \vareps Z \right) + \frac{C}{\sigma_{\min}} \frac{\vareps \log(ep)}{\sqrt{n}} \\
    & \leq C \sum_j \sup_{r \in \reals^p} \abs*{ \Exp \left[ \mathfrak{R}_{X_j, 3} \right] }
    + C \sum_j \sup_{r \in \reals^p} \abs*{ \Exp \left[ \mathfrak{R}_{Y_j, 3} \right] } 
    + \frac{C}{\sigma_{\min}} \frac{\vareps \log(ep)}{\sqrt{n}} \\
    & \leq C \nu_3 \left[
        \sum_{j=1}^{n_0}
        \frac{(\log(ep))^{3/2}}{\underline{\sigma}^3 \sqrt{n}^3}
        + \sum_{j=n_0+1}^n \left(
            \frac{1}{\vareps_j^3 p n}
            + \frac{(\log(ep))^{3/2}}{\vareps_j^3} \kappa_{j-1}(\delta_j)
        \right) 
    \right] 
    + \frac{C}{\sigma_{\min}} \frac{\vareps \log(ep)}{\sqrt{n}}. 
\end{aligned}
\end{equation*}

Based on Equation (15) in \citet{kuchibhotla2020high},
\begin{equation*}
\begin{aligned}
    \sum_{j=n_0+1}^n \frac{1}{\vareps_j^2} 
    & = \sum_{j=n_0+1}^{n} \frac{1}{\vareps^2 + \max\{\underline{\sigma}^2(n-j)-\overline{\sigma}^2,0\}} \\
    & = \sum_{i=0}^{\ceil{\overline{\sigma}/\underline{\sigma}}} \frac{1}{\vareps^2}
    + \sum_{i=\ceil{\overline{\sigma}/\underline{\sigma}}+1}^{n-n_0-1} \frac{1}{\vareps^2 + \underline{\sigma}^2 i - \overline{\sigma}^2} \\
    & \leq \frac{\ceil{\overline{\sigma} / \underline\sigma} + 1}{\vareps^2} + \frac{2}{\underline\sigma^2} \log\left(1 + \frac{\sqrt{n}\underline{\sigma}}{\vareps}\right),
\end{aligned}
\end{equation*}
\begin{equation*}
\begin{aligned}
    \sum_{j=n_0+1}^n \frac{1}{\vareps_j^3} 
    & = \sum_{i=0}^{\ceil{\overline{\sigma}/\underline{\sigma}}} \frac{1}{\vareps^3}
    + \sum_{i=\ceil{\overline{\sigma}/\underline{\sigma}}+1}^{n-n_0-1} \frac{1}{(\vareps^2 + \underline{\sigma}^2 i - \overline{\sigma}^2)^{3/2}} \\
    & \leq \frac{\ceil{\overline{\sigma} / \underline\sigma} + 1}{\vareps^3} + \frac{2}{\vareps\underline\sigma^2}.
\end{aligned}
\end{equation*}
In sum,
\begin{equation*}
\begin{aligned}
    & \mu\left( S(X_{[1,n]}), S(Y_{[1,n]}) \right) \\
    & \leq \frac{C}{\sqrt{n}} \left( \frac{\nu_3 (\log(ep))^{3/2}}{\underline{\sigma}^3} 
    + \frac{\vareps \log(ep)}{\sigma_{\min}} \right)
    + \frac{C}{pn} \nu_3 \left(\frac{\overline{\sigma}}{\vareps^3 \underline{\sigma}} + \frac{1}{\vareps\underline\sigma^2}\right) \\
    & \quad + C \nu_3 (\log(ep))^{3/2} \left( 
        \frac{\overline{\sigma}}{\vareps^2 \underline{\sigma}} + \frac{1}{\underline\sigma^2} \log\left(1 + \frac{\sqrt{n}\underline{\sigma}}{\vareps}\right)
    \right) 
    \sup_{n_0 < j \leq n} \frac{\kappa_{j-1}(\delta_j)}{\vareps_j},
\end{aligned}
\end{equation*}
which proves the lemma.

\end{appendix}

\end{document}